\documentclass[12pt,reqno]{amsart}
\usepackage{amsmath}
\usepackage{amssymb}

\thanks{
Research supported in part by
NSF Grant DMS-03-01795 (Liggett)
and by the Swedish Natural Science Research
Council, the G\"{o}ran Gustafsson Foundation (KVA)
and NSF Grant DMS-01-03841 (Steif).}

\author[Thomas M. Liggett]{Thomas M. Liggett}
\address[Thomas M. Liggett]
{Department of Mathematics,
University of California, Los Angeles
Los Angeles CA 90095, USA}
\email[Thomas M. Liggett]{tml@math.ucla.edu}

\author[Jeffrey E. Steif]{Jeffrey E. Steif}
\address[Jeffrey E. Steif]{Department of Mathematics, Chalmers University
of Technology, S-412 96 Gothenburg, Sweden}
\email[Jeffrey E. Steif]{steif@math.chalmers.se}

\usepackage{latexsym}
\usepackage{amsfonts}
\usepackage{amsmath}
\usepackage{amsthm}
\usepackage{epsfig}

\newtheorem{prop}{Proposition}[section]
\newtheorem{lemma}{Lemma}[section]
\newtheorem{thm}{Theorem}[section]
\newtheorem{cor}{Corollary}[section]
\newtheorem{definition}{Definition}[section]

\begin{document}

\title[Stochastic Domination]
{Stochastic Domination: The Contact Process,  Ising Models and FKG Measures}

\begin{abstract}
We prove for the contact process on $Z^d$, and many other
graphs, that the upper invariant measure dominates a homogeneous 
product measure
with large density if the infection rate $\lambda$ is sufficiently large.
As a consequence, this measure percolates if the corresponding
product measure percolates. We raise the question of whether domination
holds in the symmetric case for all infinite graphs of bounded degree. 
We study some asymmetric examples which we feel shed some light on this question.
We next obtain necessary and sufficient conditions for domination of a product measure for
``downward'' FKG measures. As a consequence of this general result, 
we show that the plus and minus states
for the Ising model on $Z^d$ dominate the same set of product measures.
We show that this latter fact fails completely on the homogenous 3-ary tree.
We also provide a different distinction between $Z^d$ and the homogenous 3-ary tree
concerning stochastic domination and Ising models; while it is known
that the plus states for different temperatures on $Z^d$ are never
stochastically ordered, on the homogenous 3-ary tree, almost the complete opposite is
the case. Next, we show that on $Z^d$, the set of product measures
which the plus state for the Ising model dominates is strictly increasing in the
temperature. Finally, we obtain a necessary and sufficient condition for a finite 
number of variables, which are both FKG and exchangeable, to dominate a given product measure.
\end{abstract}

\date{March 11, 2005}

\maketitle
\noindent


\section{Introduction}

There has been a significant amount of interest in determining whether
important random fields percolate for large values of some parameter.
Part of the motivation for such an interest is that some results have
been proved for parameter values above the ``percolation transition", and
one wants to make sure that such results are not vacuous. An example
is a result for the contact process in the recent paper by Broman and Steif
(2005). Another motivation involves the study of the Gibbsian nature of certain
dependent random fields -- see Maes, Redig, Shlosman and van Moffaert 
(2000), for example.

H\"aggstr\"om (1997) provides one large class of examples in which percolation
occurs. He proves that if $\mu$ is any automorphism invariant probability measure
on the bonds of the $d$-regular homogeneous tree, then $\mu$ percolates (i.e.,
there is an infinite connected component with positive probability), provided
that the marginal $\mu$-probability that an edge is present is at least
$2/d$. Such a general result does not hold on $Z^d$. Similar results
are proved for measures on the sites of the tree. For extensions to
nonamenable transitive graphs, see Benjamini, Lyons, Peres and Schramm (1999).
Another situation in which percolation has been proved for strongly
correlated random fields can be found in Bricmont, Lebowitz and Maes (1987).

One way of proving that a measure percolates is to show that it stochastically
dominates a product measure with a density that is greater than the 
critical probability for independent Bernoulli percolation. (This is not the
approach taken in the examples in the previous paragraph.) Our first
result is the following; precise definitions will be given in the relevant
sections. For any graph $G=(S,E)$, we let
$\nu_{\rho}$ denote product measure on $\{0,1\}^S$ with density $\rho$
and $\overline\nu_{\lambda}$ be the upper invariant measure 
for the contact process on $G$ with parameter $\lambda$.

\begin{thm}   \label{thm:contactdom}
Consider the graph  $Z^d$. For all $\rho <1$, there exists $\lambda$ such that
$\overline\nu_{\lambda}$ stochastically dominates $\nu_{\rho}$.
\end{thm}

\noindent{\bf Remarks.} It is easy to see that the result
for $d=1$ implies the result for $d>1$. (One sees this by comparing
the contact process on $Z^d$ with the contact process with the same
$\lambda$ on the graph that has the same vertices as $Z^d$, but only
edges in a particular coordinate direction.) This result
is the one that is relevant to Broman and Steif (2005). 
Using known results for the critical parameter for site percolation, it will follow that
$\overline \nu_{\lambda}$ percolates if $d\geq 2$ and $\lambda\geq 6.25$. 

The following question is suggested by Theorem \ref{thm:contactdom}. Consider 
the independent flip process $\zeta_t$ on $\{0,1\}^Z$ in which flips occur
from 0 to 1 at rate $r \rho$ and from 1 to 0 at rate $r(1-\rho)$, where $r>0$.
Then $\nu_{\rho}$ is invariant for $\zeta_t$. Now let
$\zeta_t$ be the stationary version of this process, chosen so that
$\zeta_t$ has distribution $\nu_{\rho}$ for all $t$ and let $\eta_t$ be the stationary contact
process with distribution $\overline\nu_{\lambda}$ at all times.
Theorem \ref{thm:contactdom} says that the two processes
can be coupled with $\zeta_t\leq \eta_t$ at a fixed time, provided that
$\lambda$ is sufficiently large.
Is it possible to construct a coupling so
that $\zeta_t\leq \eta_t$ for all times? The next proposition says that the answer
to this question is no for any choice of the parameters, except $\rho=0$. The main idea is
that a certain space-time large deviation probability is exponential
in the area of a space time box for $\zeta_t$, but exponential in
the perimeter of the box for $\eta_t$.

\begin{prop}   \label{prop:alltime}
For no parameter values except $\rho=0$ can $\{\eta_t\}$ and $\{\zeta_t\}$ be coupled so that
$P[\eta_t(x)\ge \zeta_t(x)]=1$ for all $t\ge 0$ and all $x\in Z$.
\end{prop}

\noindent{\bf Remark.} This proposition easily extends to $Z^d$.

The method of proving Theorem \ref{thm:contactdom} will allow us to prove the following 
result. We first need the following definition.

\begin{definition}
A measure $\mu$ on $\{0,1\}^S$ is {\sl downward FKG} if
\begin{equation}\label{2.2}
\text{for any finite } A\subset S, \text{the conditional measure }
\mu\{\cdot\mid
\eta\equiv 0 \text{ on } A \}\text{ is associated.}
\end{equation}
\end{definition}

In other words, if $B$ and $C$ are two increasing (resp. decreasing) subsets
of $\{0,1\}^S$, then
$$\mu\{B\cap C\mid\eta\equiv 0\text{ on }A\}\geq
\mu\{B\mid\eta\equiv 0\text{ on }A\}\mu\{
C\mid\eta\equiv 0\text{ on }A\}.$$ 

It can be shown that the FKG lattice
condition (Liggett (1985), page 78)) is equivalent to the fact that no matter how one conditions 
the configuration on $A$, the conditional distribution is associated. The word ``downward''
now refers to the fact that this is only assumed to be true when one conditions on all 
$0$'s in $A$.

\begin{thm}   \label{thm:FKG}
Let $\mu$ be a translation invariant measure on $\{0,1\}^{Z}$  which is downward FKG. 
Then the following are equivalent. \\
(1). $\mu$ dominates $\nu_{\rho}$. \\
(2). $\mu\{\eta\equiv 0\text{ on } [1,n]\}\le (1-\rho)^n$ for all $n\ge 1$.\\
(3). For all disjoint, finite  subsets $A$ and $B$  
of $\{1,2,...\}$, we have
$$
\mu\{\eta(0)=1\mid\eta\equiv 0\text{ on } A,\eta\equiv 1 \text{ on }B\}\geq \rho.
$$
\end{thm}

\noindent
{\bf Remarks.} First, we mention that it is essential here that we
are dealing with an infinite number of variables and that the process
is stationary; see also remarks after Theorem \ref{thm:exchange}. 
Next, the implications $(3) \rightarrow (1)$ 
and $(1) \rightarrow (2)$ do not require the downward FKG property. 
As far as the two corresponding reverse implications, they are both false for 
general stationary processes; for the first, see Remark 5.12 in Lyons and Steif (2003)
and for the second, see Example 1.6 in Lyons and Steif (2003).
In this latter paper, the equivalence of (1) and (3) for all conditionally negatively 
associated measures is also shown as well as (although implicitly) 
the fact that, for the family of processes studied in the Lyons-Steif paper,
domination of $\nu_\rho$ is equivalent to
$$
\mu\{\eta\equiv 1\text{ on } [1,n]\}\ge \rho^n \text{ for all } n\ge 1.
$$
In addition, Proposition \ref{prop:2sidednogood} will show that
condition (3) in this theorem cannot be modified so that
one conditions on both sides of the origin. Finally,
there is an extension of this result to $Z^d$ which we will provide. This extension
will yield the following result for the Ising model, which we think is of independent interest.

\begin{prop}   \label{prop:ising}
Fix an integer $d$ and let $\mu^{J,+}$ and $\mu^{J,-}$ be the plus and minus states
for the Ising model with nearest neighbor pair interactions on $Z^d$
with interaction parameter $J >0$.  Then for any $\rho$,
$\mu^{J,+}$ dominates $\nu_\rho$ if and only if $\mu^{J,-}$ dominates $\nu_\rho$.
\end{prop}

Interestingly, the last statement fails completely for the homogeneous 3-ary tree;
the latter graph, which we will denote by $T$, is the unique tree where every vertex has degree 3.
(Similar results will of course hold for $r$-ary trees for $r\ge 4$.)
By the free measure on $T$, we mean the Gibbs state $\mu^{J,f}$ for the Ising model
obtained by using free boundary conditions.

\begin{prop}   \label{prop:tree}
Consider the homogeneous 3-ary tree $T$ and let
 $\mu^{J,+}$, $\mu^{J,-}$ and $\mu^{J,f}$ 
be the plus, minus and free measures
for the Ising model with interaction parameter $J >0$.  If
 $\mu^{J,+}\neq \mu^{J,-}$, then there exists $0< \rho'< \rho$ such that
$\mu^{J,+}$ dominates $\nu_\rho$ but $\mu^{J,f}$ does not dominate $\nu_\rho$
and $\mu^{J,f}$ dominates $\nu_{\rho'}$ but $\mu^{J,-}$ does not dominate $\nu_{\rho'}$.
\end{prop}

We describe here another stochastic domination result for the
Ising model where the behavior is completely different depending on whether we are on $Z^d$ 
or on $T$. The first result we attribute to folklore but as we cannot find a reference, we 
include a proof here. While we have not seen this proof elsewhere, we make no claims of
its originality.

\begin{prop}   \label{prop:nodomination}
If $J_1\neq J_2$, then on $Z^d$, $\mu^{J_1,+}$ and  $\mu^{J_2,+}$ are
not stochastically ordered. 
\end{prop}

The result for $T$ is very different.

\begin{prop}   \label{prop:differentJ}
Consider the Ising model on $T$ and let $J_c$ be the critical value for $J$.
(i). If $J_c < J_1< J_2$, then
$\mu^{J_2,+}$ dominates $\mu^{J_1,+}$.  \\
(ii). For all $J_2 \ge J_c$, there exists $\alpha(J_2)$ such that
$$
\{J \in [0,J_c]: \mu^{J_2,+} \text{ dominates } \mu^{J,+} \} =[\alpha(J_2),J_c].
$$
((i) implies that $\alpha$ is a decreasing function of $J_2$).
Moreover, the smallest $J_2 > J_c$
for which $\alpha(J_2)=0$ (which corresponds to the smallest $J_2>0$ for which the plus state dominates
all plus states at lower values of $J$) is $\log(r)$ where $r$ is the unique real root of
the cubic polynomial
$$
x^3-x^2-x-1.
$$
(iii). For every $\rho<1$, there exists $J$ such that
$\mu^{+,J}$ dominates $\nu_\rho$.
\end{prop}

\noindent
{\bf Remark.} On $T$, the fact that for $J_1,J_2 \le J_c$, 
$\mu^{+,J_1}$ and $\mu^{+,J_2}$ are not stochastically ordered follows
immediately from the fact that they both have mean $1/2$
(see for example page 75 of Liggett (1985)).

Our next result concerning the Ising model and stochastic domination
compares the set of product measures that the plus state dominates at different
parameter values for $Z^d$.

\begin{prop}   \label{prop:isingdifferentJ}
For the Ising model on $Z^d$, if $0<J_1<J_2$, then
$$\sup\{\rho: \mu^{J_1,+} \text{ dominates } \nu_{\rho}\}
>
\sup\{\rho: \mu^{J_2,+} \text{ dominates } \nu_{\rho}\}.
$$
\end{prop}

\noindent
{\bf Remark.} 
Observe that Proposition \ref{prop:differentJ}(i) immediately tells us that
such a result is false on trees.

The last part of the paper gives the equivalence of 
the first two conditions of Theorem \ref{thm:FKG} 
in the context of finite exchangeable random variables which are FKG.

\begin{thm}   \label{thm:exchange}
Assume that $\eta=(\eta_1,...,\eta_n)$ is FKG and exchangeable. Then $\eta$
dominates the product measure with
density $\rho$ if and only if 
\begin{equation}\label{4.2new}
P(\eta_1=\eta_2=\cdots=\eta_n=0)\le (1-\rho)^n.
\end{equation}
\end{thm}

A large collection of examples satisfying these properties
can be obtained by taking finite pieces of an infinite exchangeable 
Bernoulli sequence $\eta_1,\eta_2,...$ (It is easy to check the FKG property
for such sequences; see Proposition 2.22 on page 83 of Liggett (1985).) 
For such a sequence, there is
a ``mixing" random variable $W$ taking values in $[0,1]$
as given in de Finetti's Theorem:
$$
W=\lim_{n\rightarrow\infty}\frac 1n(\eta_1+\cdots+\eta_n)\quad\text{a.s.}
$$
One can ask for a given $1\leq n\leq\infty$, what is 
a necessary and sufficient condition for the distribution of $(\eta_1,...,
\eta_n)$ to dominate the product measure with density $\rho$? It is easy
to see that the answer is 

(a) $EW\geq\rho$ if $n=1$,

\noindent and

(b) $W\geq \rho$ a.s. if $n=\infty$

\noindent
but Theorem \ref{thm:exchange} immediately yields, after observing that
$$
P(\eta_1=\eta_2=\cdots=\eta_n=0)=E[(1-W)^n],
$$
the following corollary.

\begin{cor}   \label{cor:exchange}
Let $\eta_1,\eta_2,...$ be an infinite exchangeable 
Bernoulli sequence with mixing random variable $W$ as above.
Then, for each $n$, the distribution of $\eta=(\eta_1,...,\eta_n)$ 
dominates the product measure with
density $\rho$ if and only if 
\begin{equation}\label{4.2}
\rho\leq 1-||1-W||_n
\end{equation}
where $||\cdot||_n$ denotes the $L_n$ norm.
\end{cor}

\noindent
{\bf Remarks.}  Note that this condition interpolates between the 
easy cases $n=1$ and $n=\infty$. 
We mention that since Theorem \ref{thm:FKG} requires an infinite number of variables,
Theorem \ref{thm:exchange} does not follow from Theorem \ref{thm:FKG}.
In addition, the natural
analogue of condition (3) in Theorem \ref{thm:FKG} is in fact not the
correct condition for stochastic domination for finite exchangeable FKG
sequences or even for sequences from
an infinite exchangeable process. For example, letting $W$ be such that $P(W=3/4)=1/2$ and
$P(W=1/4)=1/2$ and taking $n=2$, then one dominates a product measure with density
$1-(5/16)^{1/2}$ but the analogue of condition (3) only holds for densities up to 
$3/8 (< 1-(5/16)^{1/2})$. 

The rest of the paper is organized as follows. In Section \ref{sec:two},
we prove both Theorem \ref{thm:contactdom} and Proposition \ref{prop:alltime}
for the asymmetric contact process on $Z$ and also prove a result justifying
a remark after Theorem \ref{thm:FKG}. In Section \ref{sec:three}, we discuss
examples, counterexamples and extensions to more general graphs including $Z^d$.
In Section \ref{sec:four}, we prove Theorem \ref{thm:FKG} and an extension 
to $Z^d$ as well as prove Propositions \ref{prop:ising},
\ref{prop:nodomination} and \ref{prop:isingdifferentJ}.
In Section \ref{sec:five}, we prove
Propositions \ref{prop:tree} and \ref{prop:differentJ}.
In this way, the results concerning domination and Ising models 
on $Z^d$ are all in Section \ref{sec:four} while
the analogous results for trees are all in Section \ref{sec:five}.
In Section \ref{sec:six}, 
the proof of Theorem \ref{thm:exchange} is given.
In Section \ref{sec:eight}, we provide a simple example of random variables
$X_1,X_2,X_3,X_4$ which are exchangeable and FKG but do not extend to an
infinite exchangeable process. This shows that Theorem \ref{thm:exchange}
is stronger than just the statement holding for finite pieces of an
exchangeable sequence. Finally, in Section \ref{sec:nine}, we state some open questions.

\section{The One Dimensional Contact Process} \label{sec:two}
For some of the examples we have in mind, it is useful to treat first the
asymmetric contact process on $Z$. This is the continuous 
time Markov process $\eta_t$
on $\{0,1\}^Z$ in which flips at site $x$ occur from 1 to 0 at rate 1 and
from 0 to 1 at rate $\beta[p\eta(x-1)+(1-p)\eta(x+1)]$, where
$\beta>0$ and $0\leq p\leq 1$. Thus the usual
symmetric contact process with parameter $\lambda$ corresponds to
the case $\beta=2\lambda$ and $p=1/2$. Let $\overline\nu$ be the
upper invariant measure for this process (see page 135 of Liggett (1985) for this 
definition) and $\nu_{\rho}$ be the product
measure on $\{0,1\}^Z$ with density $\rho$. Recall that given two
probability measures $\mu_1, \mu_2$ on $\{0,1\}^S$,  $\mu_1$
is said to dominate $\mu_2$ (written $\mu_1\geq \mu_2$) if
$$\int hd\mu_1\geq\int hd\mu_2$$
for all increasing continuous functions $h$ on $\{0,1\}^S$. This
is equivalent to the existence of a probability measure $\gamma$ on
$\{0,1\}^S\times\{0,1\}^S$ with marginals $\mu_1,\mu_2$ that
concentrates on $\{(\eta,\zeta):\eta\geq\zeta)\}$. (See Theorem 2.4 on
page 72 of Liggett (1985).)

One way to show that $\overline\nu\geq\nu_{\rho}$  is to show that the conditional
probabilities of $\overline\nu$ satisfy 
\begin{equation}\label{2.1}
\overline\nu\{\eta(x)=1\mid\eta(y), y\neq x\}\geq \rho\quad \text{a.s.}
\end{equation}
(A generalization of this criterion is given in Liggett, Schonmann and
Stacey (1997).)
We begin by showing that (\ref{2.1}) is false if $\rho>0$, so that a
somewhat different approach is required. In this and later
arguments, we will need the following property of $\overline\nu$:
\begin{equation}\label{nu2.2}
\overline\nu \text{ satisfies the downward FKG property as defined in the introduction. }
\end{equation}
This statement is proved for the contact process on a general
graph in van den Berg,
H\"aggstr\"om and Kahn (2005). The special case in which
$B$ and $C$ is each of the form $\{\eta\equiv 0\text{ on }D\}$ was proved by
Belitsky, Ferrari, Konno and Liggett (1997). It should be noted that (\ref{2.2}) 
is not correct if the conditioning on
$\{\eta\equiv 0$ on $A\}$ is replaced by conditioning on $\{\eta(0)=1\}$. (See
Liggett (1994).) In particular, $\overline\nu$ does not satisfy the FKG lattice
condition (Liggett (1985), page 78)), since the FKG lattice condition
implies that all such conditional measures are associated.

\begin{prop}\label{prop:2sidednogood}
For any $\beta$ and $p$,
\begin{equation} \label{2.4}
\lim_{k,l\rightarrow\infty}\overline\nu\{\eta(0)=1\mid\eta\equiv 0\text{ on
} [-k,l]\backslash\{0\}\}=0.
\end{equation}
\end{prop}

\noindent
{\sl Proof:} For $k,l\geq 0$, let
$$f(k,l)=\frac {\overline\nu\{\eta(0)=1\text{ and }\eta\equiv 0\text{ on }
[-k,l]\backslash\{0\}\}}{\overline\nu\{\eta\equiv 0\text{ on }
[-k,l]\}}$$
and
$$a(l)=\overline\nu\{\eta\equiv 0\text{ on } [0,l]\}.$$
By (\ref{nu2.2}), $f(k,l)$ is decreasing
in $k$ and $l$. To see this, note for example that $f(k,l)\geq f(k,l+1)$ is
equivalent to 
\begin{eqnarray*}
\overline\nu\{\eta(0)=\eta(l+1)=0\mid\eta\equiv 0\text{ on }
[-k,l]\backslash\{0\}\}\geq\\\overline\nu\{\eta(0)=0
\mid\eta\equiv 0\text{ on }
[-k,l]\backslash\{0\}\}\overline\nu\{
\eta(l+1)=0\mid\eta\equiv 0\text{ on }
[-k,l]\backslash\{0\}\}.
\end{eqnarray*}
Since $\overline\nu$ is invariant for the process, 
$$\int Lgd\overline\nu=0$$
for all cylinder functions $g$. Here $L$ is the generator of the process.
Applying this to the indicator function of $\{\eta(y)=0\text{ for }0\leq
y\leq n\}$ and using the shift invariance of $\overline\nu$ gives
\begin{equation} \label{2.5}
\sum_{k=0}^nf(k,n-k)=\beta\big[pf(0,n+1)+(1-p)f(n+1,0)\big]\frac{a(n+1)}
{a(n)}\leq\beta.
\end{equation}
By the monotonicity of $f$, the limit in 
(\ref{2.4}) exists. Call it $\gamma$.
By (\ref{2.5}),
$$\gamma(n+1)\leq\beta.$$
It follows that $\gamma=0$ as required. $\blacksquare$

Next, we show that the situation is quite different if one conditions
on the configuration on one side of the origin rather than on both sides.
This will enable us to prove domination of a product measure by
constructing the coupling measure sequentially. Another class of
processes for which domination of nontrivial product measures is proved via bounds
on one sided conditional probabilities, while (\ref{2.1}) fails, are the measures $P^f$
studied by Lyons and Steif (2003), in which the function $f$ has a positive
geometric mean, but 0 harmonic mean.

\begin{prop} \label{prop2.6}
For all $\beta\geq 4, 0\leq p\leq 1$, and all $l\geq 0$,
\begin{equation} \label{2.7}
\overline\nu\{\eta(0)=1\mid\eta\equiv 0\text{ on } [1,l]\}\geq \frac{\beta-4}{\beta}.
\end{equation}
\end{prop}

\noindent
{\sl Proof:} The idea of the proof is to use the known fact that
the probability of having $n$ consecutive 0's is exponentially
small with some definite rate. To translate this into a statement about
conditional probabilities, we need some monotonicity.
Let $f(l)=\overline\nu\{\eta(0)=0\mid\eta\equiv 0\text{ on
} [1,l]\}$. Now $f(l)$ is increasing in $l$ by (\ref{nu2.2}). Note that
\begin{equation} \label{2.8}
\prod_{l=0}^nf(l)=a(n)
\end{equation}
where $a(n)$ is defined as in the previous proposition. For $\beta\geq 4$, let
$\mu$ be the stationary renewal measure on $\{0,1\}^Z$ defined by
saying that the numbers of 0's between successive 1's are i.i.d. random
variables $\xi_i$ with tail probabilities given by
\begin{equation} \label{2.9}
P(\xi_i\geq k)=\frac{(2k)!}{k!(k+1)!}\beta^{-k},\quad k\geq 0.
\end{equation}
Holley and Liggett (1978) proved 
\begin{equation} \label{2.10}
a(n)\leq\mu\{\eta\equiv 0\text{ on }[0,n]\},\quad n\geq 0.
\end{equation}
(This statement for $n=0$ is (2.1) in their paper; for general $n$
it comes from (2.2) applied to the set $A=\{0,...,n\}$.)
Combining (\ref{2.8}), (\ref{2.10}) and the monotonicity of $f$, we see that
\begin{equation} \label{2.11}
\lim_{n\rightarrow\infty}f(n)=\lim_{n\rightarrow\infty}[a(n)]^{1/n}
\leq \liminf_{n\rightarrow\infty}[\mu\{\eta\equiv 0\text{ on }[0,n]\}]^{1/n}.
\end{equation}
By (\ref{2.9}),
$$\mu\{\eta\equiv 0\text{ on }[0,n]\}=\mu\{\eta(0)=1\}\sum_{k=n+1}^{\infty}
\frac{(2k)!}{k!(k+1)!}\beta^{-k}.$$
This, together with (\ref{2.11}) and Stirling's formula, implies
$$\lim_{n\rightarrow\infty}f(n)\leq 4/\beta.$$
This proves (\ref{2.7}), since $f(n)$ is increasing. $\blacksquare$

\noindent{\bf Remarks.} Exponential decay of $a(n)$ has been proved for
symmetric contact processes on $Z^d$ for all $\lambda$ strictly above the 
critical value. (See (1.8) on page 36 and Theorem 2.30(b) on page 57
of Liggett (1999).) (For $Z^d$, $a(n)$ refers to the probability that there are 0's at
$n$ specified locations; while this probability depends on the $n$ locations,
the exponential rate is uniform over all such sets.)
The advantage of using (\ref{2.10}) instead is that
it works for asymmetric processes as well, and that it gives an explicit decay rate.
\medskip

Now, we need to show that the conditioning in (\ref{2.7}) on having all 0's 
at sites adjacent to the origin is the worst case.

\begin{prop} \label{prop:2.12}
Suppose $A$ and $B$ are disjoint, finite
subsets of $\{1,2,...\}$. Then
\begin{equation} \label{2.13}
\overline\nu\{\eta(0)=1\mid\eta\equiv 0\text{ on
} A,\eta\equiv 1 \text{ on }B\}\geq \frac{\beta-4}{\beta}.
\end{equation}
\end{prop}

\noindent
{\sl Proof:} By (\ref{nu2.2}),
\begin{eqnarray*} 
\overline\nu\{\eta(0)=1,\eta\equiv 1\text{ on }B\mid\eta\equiv 0\text{ on
} A\}
\geq \\\overline\nu\{\eta(0)=1\mid\eta\equiv 0\text{ on }
A\}\overline\nu\{\eta\equiv 1\text{ on }B\mid\eta\equiv 0\text{ on }
A\}
\end{eqnarray*} 
so $$\overline\nu\{\eta(0)=1\mid\eta\equiv 0\text{ on
} A,\eta\equiv 1 \text{ on }B\}\geq\overline\nu\{\eta(0)=1\mid\eta\equiv 0\text{ on
} A\}.$$
Using (\ref{nu2.2}) again, we see that
$$\overline\nu\{\eta(0)=1\mid\eta\equiv 0\text{ on
} A\}\geq \overline\nu\{\eta(0)=1\mid\eta\equiv 0\text{ on
} [1,l]\}$$
for any $l$ such that $A\subset [1,l]$. So, (\ref{2.13}) follows from 
(\ref{2.7}).
$\blacksquare$

The following statement follows immediately from Proposition \ref{prop:2.12}, since
one can construct the required coupling measure sequentially.
\begin{thm}   \label{thm:2.14}
If $\beta\geq 4$ and $0\leq p\leq 1$, then
$$\overline\nu\geq\nu_{\rho}$$
for $\rho=(\beta-4)/\beta.$
\end{thm}

\noindent{\bf Remark.} In Theorem \ref{thm:2.14}, one can equally well consider the
stationary distribution for oriented percolation, which is a discrete
time version of the contact process. The main ingredients of the proof
of the theorem in continuous time are (\ref{nu2.2}) and (\ref{2.10}). 
Van den Berg,
H\"aggstr\"om and Kahn (2005) is mainly devoted to the discrete time
setting, and the continuous time results are deduced from them. The
discrete time analogue of (\ref{2.10}) is given in Liggett (1995).
\medskip

We now turn to the \\
\noindent
{\sl Proof of Proposition \ref{prop:alltime}:}

If the coupling  $\zeta_t\leq\eta_t$ were possible for all times, then for all $N,T$,
\begin{equation}\label{2.15}
P(\zeta_t(n)=0\text{ for all }1\leq n\leq N\text{ and all }0\leq t\leq T)\geq 
\end{equation}
\begin{equation}\label{2.16}
P(\eta_t(n)=0\text{ for all }1\leq n\leq N\text{ and all }0\leq t\leq T).
\end{equation}

To see that this is not possible, first compute (\ref{2.15}):
$$P(\zeta_t(n)=0\text{ for all }1\leq n\leq N\text{ and all }0\leq t\leq T)
=(1-\rho)^Ne^{-r\rho T N}.$$
To see this, note that the event in question occurs if and only if the
configuration is $\equiv 0$ on $\{1,..., N\}$  at time 0 and none of these
sites flips to a 1 in time $T$. Note that the right side above is
exponentially small in the ``area" $NT$ of the space--time rectangle
$\{1,...,N\}\times [0,T]$. 

It remains to check that (\ref{2.16}) is bounded below by a quantity that 
is exponentially small in the perimeter $N+T$ of this rectangle. In doing
so, we will use the graphical representation of the contact process
-- see pages 172--174 of Liggett (1985) or pages 32--33 of Liggett (1999)
for its description. Letting $A_{N,T}$ denote the event that
there are no arrows in the graphical representation
from 0 to 1 or from $N+1$ to $N$ during the time period $[0,T]$,
(\ref{2.16}) is bounded below by
\begin{equation}
P(\{\eta_0(n)=0\text{ for }1\leq n\leq N\}\cap A_{N,T}).
\end{equation}
This, in turn, is bounded below by
$$P(\eta_0(n)=0\text{ for }1\leq n\leq N)e^{-2\lambda T}
\geq[\overline\nu\{\eta:\eta(0)=0\}]^Ne^{-2\lambda T},$$
where the inequality comes from the unconditioned version of (\ref{nu2.2}).
Taking $N=T$, we see that (\ref{2.15}), (\ref{2.16}) cannot hold for all $N,T$.
$\blacksquare$

\section{Examples, Extensions, and Counterexamples}\label{sec:three}
One can deduce from Theorem \ref{thm:2.14} domination of product measures, and hence
percolation, in many situations when the corresponding product measures percolate.
Here are some examples:

(a) The upper invariant measure $\overline\nu$ for the symmetric contact process on
$Z^d$ with parameter $\lambda$ dominates $\nu_{\rho}$ with
$\rho=(\lambda-2)/\lambda$ if $\lambda\geq 2$. (This gives us Theorem \ref{thm:contactdom}.)
Therefore, using the
upper bound of $.679492$ for the critical value of site percolation on $Z^2$
(Wierman (1995)), we conclude that $\overline\nu$ percolates if
$\lambda\geq 6.25$ and $d\geq 2$. The domination statement for $d=1$ is just
Theorem \ref{thm:2.14} with $p=1/2$ and $\beta=2\lambda$. 

(b) The upper invariant measure $\overline\nu$ for the symmetric contact
process on $\{0,1,...\}$ with parameter $\lambda$ dominates $\nu_{\rho}$
with $\rho=(\lambda-4)/\lambda$ if $\lambda\geq 4$. To see this, compare
this process with the asymmetric contact process on $Z$ with $\beta=\lambda$
and $p=0$.

(c) Let $G$ be any graph with the property that after deleting a set
of edges (but no vertices), the resulting graph is a union of 
disconnected copies of $\{0,1,...\}$. Then
the upper invariant measure $\overline\nu$ for the symmetric contact
process on $G$ with parameter $\lambda$ dominates $\nu_{\rho}$ with
$\rho=(\lambda-4)/\lambda$. An example is the graph in which $n$
semi-infinite spokes meet at a common vertex. Another is any infinite
tree with no leaves. Note that in these examples, the
value of $\rho$ does not depend on the complexity of the graph.

\medskip
\noindent {\bf Remark.} Of course, percolation for large $\lambda$ on 
homogeneous trees follows from H\"aggstr\"om's 1997 result. Example
(c) in the case of homogeneous trees gives domination of a product
measure of large density as well.

\medskip
Based on example (c) above, one might think that the stated domination 
holds for any graph that contains a copy of $\{0,1,...\}$. We cannot show this
to be the case for the symmetric contact process; see Section \ref{sec:six}.
The remainder of this section
is devoted to studying certain examples of asymmetric contact processes;
we feel that these examples shed some light on what might be the difficulties in
proving the above.

Let $G$ be the graph consisting of a copy of $\{0,1,...\}$, together with
$n$ other vertices $v_1,...,v_n$  that have edges only to $0$. The infection rate
at vertex $k$ is $\lambda\eta(k+1)$ for vertices in $\{0,1,...\}$, and is
$\lambda\eta(0)$ for $v_1,...,v_n$.  Let $\overline\nu$ be the upper invariant
measure for the process. Note that $\overline\nu$ restricted to
$\{x_k,0,1,...\}$ dominates $\nu_{\rho}$ for $\rho=(\lambda-4)/\lambda$
by example (b) above. However, for any $\lambda$, the largest value of $\rho$
so that $\overline\nu\geq\nu_{\rho}$ tends to zero as $n\rightarrow\infty$:

\begin{prop}   \label{prop:3.1}
Let
$$\rho(\lambda,n)=\sup\{\rho:\overline\nu\geq\nu_{\rho}\}.$$
Then, with $G$ as in the previous paragraph,
$$\text{ (a) }\,\,\,\,\,\, \lim_{n\rightarrow\infty}\rho(\lambda,n)=0$$
for each $\lambda$, and
$$
\text{ (b) }\,\,\,\,\,\,\lim_{\lambda\rightarrow\infty}\rho(\lambda,n)=1$$
for each $n$.
\end{prop}

\noindent
{\sl Proof:}
To prove part (a), fix $\lambda$ and $\rho$ and choose $M$ such that
$$
e^{-M}=\rho/2.
$$
Let $A$ be the event that $\eta_t(0)=0$ for all $t\in [0,M]$,
$B$ be the event that $\eta_M(v_k)=0$ for $k=1,\ldots,n$ and
$C$ be the event that in the graphical representation,
for all $v_k$, $k=1,\ldots,n$, there is a death during $[0,M]$.
Observe that $P(A)$ is independent
of $n$ because of the 1-way nature of the infection. Now
$$
P(B)\ge P(A\cap B)\ge P(A\cap C)
$$
$$
=P(A)(1-e^{-M})^n=P(A)(1-[\rho/2])^n,
$$
since $A$ and $C$ are independent.
If $\overline\nu$ did dominate $\nu_{\rho}$, then we would have
$(1-\rho)^n\ge P(B)$, which gives
$$
(1-\rho)^n\ge P(A)(1-[\rho/2])^n.
$$
Since $P(A)$ does not depend on $n$, the above fails for large $n$.

For part (b), we will show that $\overline\nu\geq\nu_{\rho}$ if
\begin{equation}\label{3.2}
4n\leq\lambda (1-\rho)^n.
\end{equation}
To do so, take $\eta$ to be $\overline\nu$ distributed
and $\zeta$ to be $\nu_{\rho}$ distributed. By Proposition
\ref{prop:2.12},
$$P(\eta(v_k)=1\mid \eta(j),j\geq 0)\geq \frac {\lambda-4}{\lambda}\quad
\text{a.s.}$$
for each $k=1,...,n$. Therefore,
\begin{equation}\label{3.3}
\aligned P(\eta(v_k)=1\ &\forall\ 1\leq k\leq n\mid \eta(j), j\geq
0)\geq1-(4n/\lambda)\\
&\geq1-(1-\rho)^n=P(\zeta(i)=1
\text{ for some }1\leq i\leq n\}\quad\text{a.s.,}\endaligned
\end{equation}
where the second inequality comes from (\ref{3.2}). By Theorem \ref{thm:2.14} and 
(\ref{3.2}),
the distribution of $(\eta(j),j\geq 0)$ dominates $\nu_{\rho}$.
Combining this with (\ref{3.3}) gives $\overline\nu\geq\nu_{\rho}$ as
required. $\blacksquare$

We next let $G$ be the graph $\{0,1,...\}$, to which are added $n$ neighbors to
vertex $n$ for each $n\geq 1$. The infection rate is $\lambda\eta(n+1)$
at vertex $n$ and is $\lambda\eta(n)$ at the $n$ neighbors that were
added to vertex $n$.  

\begin{cor}   \label{cor:3.4}
For the graph $G$ in the last paragraph,
$\overline\nu$ does not dominate $\nu_{\rho}$ for any $\rho>0$.
\end{cor}
\noindent
{\sl Proof:} 
The result follows from part (a) of Proposition \ref{prop:3.1},
since for each $n$, the process restricted to $\{n,n+1,...\}$,
together with the $n$ vertices that were attached to $n$, is a copy of the 
process considered in Proposition \ref{prop:3.1}. $\blacksquare$

\section{Domination characterization and the Ising model}\label{sec:four}

Since the extension of Theorem \ref{thm:FKG} to $Z^d$
is slightly messy, we choose to first prove this result
for $Z$ and afterwards state the result for $Z^d$ and 
outline the proof.

\noindent
{\sl Proof of Theorem \ref{thm:FKG}:}  \\
(1) implies (2) is trivial. For (2) implying (3), (2) says that for all $n$
$$
 \mu\{\eta(0)=0\}\prod_{i=1}^{n-1} \mu\{\eta(0)=0\mid\eta\equiv 0\text{ on} [1,i]\}
\le (1-\rho)^n.
$$
The assumption of downward FKG easily gives that
$\mu\{\eta(0)=0\mid\eta\equiv 0\text{ on } [1,i]\}$ is increasing in $i$ and hence
converges to a limit $L$. By the above, we must have $L\le 1-\rho$. It follows that
 $\mu\{\eta(0)=0\mid\eta\equiv 0\text{ on} [1,i]\}$ is $\le 1-\rho$ for each $i$ and
(3) follows as in the proof of Proposition \ref{prop:2.12} using the
downward FKG assumption. Finally, (3) implies (1) by constructing the coupling
measure sequentially as in Theorem \ref{thm:2.14}. $\blacksquare$

The extension to general $d$ requires the notion of a lexicographical order. For simplicity,
we do this only for $d=2$; the reader can easily extend to general $d$. 
Let $\mathcal P$ (for past) be
the subset of $Z^2$ given by $\{(i,j): \{i <0 \text{ and } j\le 0\} \text{ or }j < 0\}$. 
(This is the set of vertices
below the $x$-axis or to the left of 0 on the $x$-axis.) For $x\in Z^2$, 
let ${\mathcal P}_x={\mathcal P}+x$. 

\begin{thm}   \label{thm:FKGd}
Let $\mu$ be a translation invariant measure on $\{0,1\}^{Z^2}$  which is downward FKG. 
Then the following are equivalent. \\
(1). $\mu$ dominates $\nu_\rho$. \\
(2). $\mu\{\eta\equiv 0\text{ on } [1,n]^2\}\le (1-\rho)^{n^2}$ for all $n$.\\
(3). For all disjoint, finite  subsets $A$ and $B$ of ${\mathcal P}$, we have
\begin{equation} 
\mu\{\eta((0,0))=1\mid\eta\equiv 0\text{ on } A,\eta\equiv 1 \text{ on }B\}\geq \rho.
\end{equation}
\end{thm}

\noindent
{\sl Proof:} \\
(1) implies (2) is trivial. (3) implies (1) as in the previous proof after one observes that
on a square box, one can order the vertices $(a_1,a_2,\dots)$ in such a way that
for all $i$, $\{a_1,a_2,\dots,a_{i-1}\}\subset P_{a_i}$; this allows one to do the coupling
sequentially. (The ordering is of course just going from left to right starting on the bottom row
and working upwards.) For (2) implies (3), we proceed as follows.
For any $n\ge 1$ and $x\in  [1,n]^2$, let $A^n_x= [1,n]^2\cap {\mathcal P}_x$.
From the order used above, (2) immediately implies that for all $n\ge 1$
\begin{equation} \label{eqn:product}
\prod_{v\in [1,n]^2} \mu\{\eta(v)=0\mid\eta\equiv 0\text{ on } A^n_v\}
\le (1-\rho)^{n^2}.
\end{equation}
Let 
$$
L=\lim_{n\to \infty}{\mu\{\eta((0,0))=0\mid\eta\equiv 0\text{ on } {\mathcal P}\cap [-n,n]^2\}}.
$$
The downward FKG condition easily implies, as before, that
(a) the terms in the above limit are nondecreasing in $n$
and hence the limit $L$ exists and (b)
for all disjoint, finite  subsets $A$ and $B$ of ${\mathcal P}$, we have
\begin{equation} 
\mu\{\eta((0,0))=0\mid\eta\equiv 0\text{ on } A,\eta\equiv 1 \text{ on }B\}\le L.
\end{equation}
(To see, for example, (a), one takes $A={\mathcal P}\cap [-n,n]^2$,
$B=\{(0,0)\}$ and  $C={\mathcal P}\cap [-n-1,n+1]^2\backslash [-n,n]^2$.)
(3) will therefore be proved if we show that $L\le 1-\rho$. To show this, first note that, again by 
downward FKG, there is a uniform lower bound $b$ on all the factors appearing in the left 
side of (\ref{eqn:product}). (Of course $b=\mu\{\eta((0,0))=0\}$.) If  $L> 1-\rho$, choose $\delta >0$ 
so that $L-\delta > 1-\rho$ and then choose $N$ so that
$$
\mu\{\eta((0,0))=0\mid\eta\equiv 0\text{ on } {\mathcal P}\cap [-N,N]^2\}\ge L-\delta.
$$
Choose $r\in (0,1)$ so that
$$
b^{1-r}(L-\delta)^r > 1-\rho
$$
and finally choose $M$ so large that that number of $x$'s in $[1,M]^2$ for which 
$$
x+[{\mathcal P}\cap [-N,N]^2]\subseteq A^M_x
$$
is at least $rM^2$. It is clear geometrically that this can be done. By choice of $M$,
it follows, again from the downward FKG condition, that
$$
\prod_{v\in [1,M]^2} \mu\{\eta(v)=0\mid\eta\equiv 0\text{ on } A^M_v\}
\ge b^{(1-r)M^2}(L-\delta)^{rM^2}.
$$
This is strictly larger than $(1-\rho)^{M^2}$ contradicting (\ref{eqn:product}) for $n=M$.
$\blacksquare$

\noindent
{\bf Remarks.} (1). One can see from the proof that the above conditions are equivalent to requiring 
that (2) holds only for all sufficiently large $n$. \\
(2). Theorem \ref{thm:FKGd} together with (\ref{nu2.2}) and the remark after Proposition 
\ref{prop2.6} concerning exponential decay in the symmetric contact process 
immediately yields the following corollary.

\begin{cor}   \label{cor:french}
For the symmetric contact process on $Z^d$, for any $d$, and any 
$\lambda > \lambda_c$, we have that $\overline\nu_{\lambda}$ 
stochastically dominates $\nu_{\rho}$ for some $\rho >0$.
\end{cor} 

We have decided, in order to save space, not to define the Ising model 
which perhaps most readers are familiar with; for those who are not,
see Liggett (1985), Chapter 4 for all definitions and claims. To stick with 0,1 valued
random variables, we are representing the state -1 in the Ising model by 0. We deal exclusively with
the Ising model with no external field and hence the only parameter (besides the graph on which the model
lives) is the interaction parameter $J$. We will write $\mu^{J,+}$ and $\mu^{J,-}$ for
the plus and minus states at interaction $J$ and we will let $J_c$ denote the critical parameter.
Here the underlying graph could be $Z^d$, the homogeneous 3-ary tree $T$ or any other graph.

We now prove Proposition \ref{prop:ising}. 

\noindent
{\sl Proof of Proposition \ref{prop:ising}:}  \\
We prove this only for $d=2$, the proof for general $d$ being identical.
It is known that $\mu^{J,+}$ and $\mu^{J,-}$ both satisfy the FKG condition. 
It is also immediate that both of these measures satisfy the uniform finite energy property
which means that for some $\gamma>0$, the probability of having a 1 or having a 0 at
a site conditioned on everything else is always at least $\gamma$.

Assume now that $\mu^{J,+}$ dominates $\nu_\rho$. Fix $\epsilon >0$. We have that
$$
\mu^{J,+}\{\eta\equiv 0\text{ on } [1,n]^2\}\le (1-\rho)^{n^2}
$$
Let $B_n= [0,n+1]^2 \backslash [1,n]^2$. We have
$$
\mu^{J,+}\{\eta\equiv 0\text{ on } [1,n]^2\} \ge
\mu^{J,+}\{\eta\equiv 0\text{ on } [0,n+1]^2\}=
$$
$$
\mu^{J,+}\{\eta\equiv 0\text{ on } B_n\}
\mu^{J,+}\{\eta\equiv 0\text{ on } [1,n]^2\mid\eta\equiv 0\text{ on } B_n\} \ge
$$
$$
\gamma^{|B_n|}
\mu^{J,-}\{\eta\equiv 0\text{ on } [1,n]^2\mid\eta\equiv 0\text{ on } B_n\} \ge
$$
$$
\gamma^{|B_n|}\mu^{J,-}\{\eta\equiv 0\text{ on } [1,n]^2\}.
$$
All of these equalities and inequalities follow immediately from well known results 
about the Ising model. For example, for the second to last inequality,
the second factors are equal because the plus and minus states are Markov 
fields for the same set of conditional probabilities and the first factors
satisfy this inequality by the definition of $\gamma$. This gives
$$
\mu^{J,-}\{\eta\equiv 0\text{ on } [1,n]^2\}
\le \gamma^{-|B_n|} \mu^{J,+}\{\eta\equiv 0\text{ on } [1,n]^2\} 
$$
which is at most $\gamma^{-|B_n|}(1-\rho)^{n^2}$. For large $n$, this is at most
$[(1-\rho)(1+\epsilon)]^{n^2}$. By Theorem \ref{thm:FKGd} and Remark (1) after the proof, we conclude
that $\mu^{J,-}$ dominates $\nu_{w}$ where $w=1-[(1-\rho)(1+\epsilon)]$. As $\epsilon$ is arbitrary,
we are done. 
$\blacksquare$

We finally point out that condition (3) in Theorem \ref{thm:FKGd} can also be used to obtain some
upper bounds on those $\rho$ for which $\nu_\rho$ is dominated by $\mu^{J,h}$ where $\mu^{J,h}$ 
is the Ising model with interaction parameter $J$ and external field $h$. To do this, one could for
example, for $d=2$, place $0$'s at locations $\{(-1,0),(0,-1),(1,-1)\}$ and get some upper bounds
on the conditional probability that the origin is in state 1. By Theorem \ref{thm:FKGd},
this would give bounds on what product measure densities one would dominate. We illustrated this 
method in 2 dimensions but this of course could be done
in any dimension. However, in 2 dimensions, because of the exact formula for the pressure 
for the Ising model due to Onsager, we can write down explicitly the optimal $\rho=\rho(J)$.

\begin{cor}   \label{cor:onsager}
For the 2 dimensional Ising model with parameter $J$, the maximal $\rho$ for which $\mu^{+,J}$ 
dominates $\nu_\rho$ is
$$
1-\frac{e^{2J}}{2e^{\sigma(J)}}
$$
where 
$$
\sigma(J)=\frac{1}{2\pi^2} \int_0^\pi\int_0^\pi \log[\cosh^2 2J-\sinh 2J(\cos x +\cos y)]dxdy.
$$
\end{cor}

\noindent
{\sl Proof:} 
For this proof, it is easier to deal with $\pm 1$ variables.
Let 
$$
Z_n(J)=\sum_{\eta\in \{\pm 1\}^{[1,n]^2}} e^{J U_1 -J U_2}
$$
where $U_1$ is the number of unordered nearest neighbor pairs in
$[1,n]^2$ where $\eta$ agrees and
where $U_2$ is the number of unordered nearest neighbor pairs in
$[1,n]^2$ where $\eta$ disagrees. This is just the usual normalization
(partition function) for the Ising model on
$[1,n]^2$. Onsager's formula (see Thompson (1972), page 132))
says that 
$$
\lim_{n\to\infty}(Z_n(J))^{1/n^2}=2e^{\sigma(J)}
$$
with $\sigma(J)$ as above. Next, if $\mu^{J,+,n}$ is the Ising model on 
$[1,n]^d$ with parameter $J$ and plus boundary conditions, then
$$
\mu^{J,+,n}(\{\eta\equiv 0\text{ on } [1,n]^2\})=\frac{e^{JU_3-JU_4}}{Z^+_n(J)}
$$
where $U_3$ is the number of unordered nearest neighbor pairs in
$[1,n]^2$, $U_4$ is the number of unordered nearest neighbor pairs with
exactly one point is in $[1,n]^2$ and $Z^+_n(J)$ is the normalization
needed to make $\mu^{J,+,n}$ a probability measure.
Since $U_3=2n^2 +o(n^2)$, $U_4=O(n)$ and
$$
\lim_{n\to\infty}(Z_n(J))^{1/n^2}=
\lim_{n\to\infty}(Z^+_n(J))^{1/n^2}
$$
(see Georgii (1988), pages 322) we have that
$$
\lim_{n\to\infty} 
\mu^{J,+,n}(\{\eta\equiv 0\text{ on } [1,n]^2\})^{1/n^2}=\frac{e^{2J}}{2e^{\sigma(J)}}.
$$
It is elementary to check that this last statement is equivalent to
$$
\lim_{n\to\infty} 
\mu^{J,+}(\{\eta\equiv 0\text{ on } [1,n]^2\})^{1/n^2}=\frac{e^{2J}}{2e^{\sigma(J)}}.
$$
It next easily follows from Theorem \ref{thm:FKGd} together with the remark 
immediately afterwards that $\mu^{+,J}$ dominates $\nu_\rho$ for
$$
\rho =1-\frac{e^{2J}}{2e^{\sigma(J)}}
$$
but for no larger $\rho$.
$\blacksquare$

\noindent{\bf Remark.}
Of course the above proof shows that for the Ising model with interaction 
parameter $J$ and external field $h$ on $Z^d$, the maximal $\rho$ for which the measure 
dominates $\nu_\rho$ is
$$
1-\frac{e^{dJ}}{P(J,h)}
$$
where $P(J,h)$ is the limiting $n^d$th root of the partition function on the box 
$[1,n]^d$. It is just that we only have a formula for $P(J,h)$ when $d=2$ and $h=0$.

\medskip
\noindent
{\sl Proof of Proposition \ref{prop:nodomination}:}  \\
For this proof, it is again simplest to think of the model as $\pm 1$ valued.
Assume that $J_1 < J_2$ and let $m_1$ and $m_2$ be the expected values of a fixed spin
under $\mu^{+,J_1}$ and $\mu^{+,J_2}$. It is known that $m_1\le m_2$
(see page 186 of Liggett (1985)). If $m_1=m_2$,
then it is clear that they cannot be ordered (see page 75 of Liggett (1985)).
If $m_1< m_2$, choose $m^*\in (m_1,m_2)$
and consider the event 
$$
E_n=\{\sum_{x\in [-n,n]^d} \eta(x) < -m^* (2n+1)^d\}.
$$
Known results (see Schonmann (1987) where the reader is reminded of
the first below inequality and where the second inequality is proved) for large 
deviations for the Ising model tell us that
$$
\mu^{+,J_1}(E_n) \le c_1 e^{-c_2n^d}
$$
and
$$
\mu^{+,J_2} (E_n)\ge c_3 e^{-c_4n^{d-1}}
$$
for strictly positive constants $c_1,c_2,c_3$ and $c_4$. 
Taking $n$ large, one sees that
$\mu^{+,J_1}$ is not dominated by $\mu^{+,J_2}$. 
$\blacksquare$

\noindent
{\sl Proof of Proposition \ref{prop:isingdifferentJ}:}  \\
By Theorem 4.1 and Remark 1 following its proof, the set of $\rho$'s for
which
$\mu^{J,+}\geq
\nu_{\rho}$ is determined by
$$\limsup_{n\rightarrow\infty}\mu^{J,+}\{\eta\equiv 0\text{ on }[1,n]^d\}^{1/n^d}.$$
(Although we do not need it, this $\limsup$ is in fact a limit by an
easy subadditivity argument using the FKG inequality.)
Therefore, the assertion in the proposition is just the statement that this
$\limsup$ is strictly increasing in $J$. From the proof of Proposition 1.2,
we see that this $\limsup$ is the same as
$$\limsup_{n\rightarrow\infty}\mu^{J,-}\{\eta\equiv 0\text{ on }[1,n]^d\}^{1/n^d},$$
and then by interchanging the roles of zeros and ones, it is the same as
\begin{equation}\label{final}
\limsup_{n\rightarrow\infty}\mu^{J,+}\{\eta\equiv 1\text{ on }[1,n]^d\}^{1/n^d}.
\end{equation}

We will deduce the strict monotonicity of this quantity from Griffiths' inequality. (See
page 186 of Liggett (1985), for example). To do so, let
$$\chi_A(\eta)=\prod_{x\in A}[2\eta(x)-1]$$
for finite $A\subset Z^d$ and $\eta\in\{0,1\}^{Z^d}.$ These form an orthonormal basis
for $L_2(\nu_{\frac 12}$). For a finite $C\subset Z^d$, expand $f_C=
1_{\{\eta:\eta\equiv 1\text {on }C\}}$ in this basis:
$$f_C(\eta)=\frac 1{2^{|C|}}\sum_{A\subset C}\chi_A(\eta).$$
Now we use Griffiths' inequality and its formalism. Let $\Lambda$
be a large box in $Z^d$ that contains $C$, and let $\mu^{J,\Lambda}$ be the Gibbs
state on $\Lambda$ with plus boundary conditions.
The sums on $B$ below are over nearest
neighbor pairs contained in $\Lambda$ and singletons in $\Lambda$
with a neighbor outside $\Lambda$. We denote the
covariance with respect to $\mu^{J,\Lambda}$ by cov$_J$. Then, using Griffiths'
inequality in the final step, we see that
\begin{eqnarray*}
\frac{d}{dJ}\int f_C d\mu^{J,\Lambda}=\frac 1{2^{|C|}}\sum_{A\subset C}
\frac{d}{dJ}\int \chi_Ad\mu^{J,\Lambda}= 
\frac 1{2^{|C|}}\sum_{A\subset C}
\sum_B \text{ cov}_J(\chi_A,\chi_B)
\\=\sum_B\text{ cov}_J(f_C,\chi_B)
\geq \int f_Cd\mu^{J,\Lambda}\sum_{B\subset C}\bigg[1-\int\chi_Bd\mu^{J,\Lambda}\bigg].
\end{eqnarray*}
Dividing this inequality by $\int f_C d\mu^{J,\Lambda}$, integrating with respect to
$J$, and then passing to the limit as $\Lambda$ exhausts $Z^d$, we see that
for $0<J_1<J_2$,
\begin{eqnarray*}
\log\int f_Cd\mu^{J_2,+}-\log\int f_Cd\mu^{J_1,+} &\geq &   
4(\# \text{ nearest neighbor pairs in }C)  \\
&&\int_{J_1}^{J_2}\mu^{J,+} \{\eta:\eta(0)=1,\eta(e)=0\}dJ, 
\end{eqnarray*}
where $e$ is a neighbor of 0. Applying this to $C=[1,n]^d$ gives the
strict monotonicity of (\ref{final}) as required.
$\blacksquare$

\section{Ising models on trees and domination} \label{sec:five}

Throughout this section, $T$ will denote the
homogeneous 3-ary tree whose vertices are $V(T)$ and edges are
$E(T)$. We first need to define a 2 state {\it tree indexed Markov chain}. 
Let $\{P(i,j)\}_{i,j\in \{0,1\}}$ be the transition matrix for an irreducible
2 state Markov chain with stationary distribution $\pi$.
From this we will define a probability measure $\mu$ on $\{0,1\}^{V(T)}$. 
Fixing a connected set $A\subseteq V(T)$ and an $\eta \in \{0,1\}^{A}$, we define $\mu(\eta)$
as follows. Choose an arbitrary
element $a\in A$. Let $F$ be the set of directed edges $(x,y)$ where $x,y\in A$ and $x$ is closer
to $a$ than $y$ is. Now, define
$$
\mu(\eta)= \pi(\eta(a))\prod_{(x,y)\in F} P(\eta(x),\eta(y)).
$$
Using the fact that any 2 state Markov chain is reversible, it is easy to check that this
definition is independent of the choice of $a$ and also that $\mu$ defined for different $A$'s
as above are consistent. It is also easy to check that (i) on
any biinfinite line through the tree, we see a copy of the above stationary Markov chain
and that (ii) $\mu$ is invariant under all tree automorphisms.

Before giving the proof of Proposition \ref{prop:tree}, we need the following
result analogous to Theorem \ref{thm:FKGd}. Fix an origin $o\in V(T)$ and 
let $T_n$ be the induced subtree of $T$ whose vertices are the elements in $V(T)$
within distance $n$ of $o$. The vertices of $T_n$ will be denoted by $V(T_n)$.

\begin{prop}   \label{prop:treeFKG}
Let $\{P(i,j)\}_{i,j\in \{0,1\}}$ be a transition matrix with $P(0,1)\le P(1,1)$
(or equivalently $P(1,0)\le P(0,0)$) and let $\mu$ be the distribution of the
corresponding tree indexed process. Then the following are equivalent. \\
(1). $\mu$ dominates $\nu_\rho$. \\
(2). $\mu\{\eta\equiv 0\text{ on } V(T_n)\}\le (1-\rho)^{|V(T_n)|}$ for all $n$.\\
(3). $P(0,1)\ge \rho$.
\end{prop}

\noindent
{\sl Proof:}
(1) implies (2) is trivial. For (3) implies (1), note that our assumption that
$P(0,1)\le P(1,1)$ implies $\pi(1)\ge \rho$.
Fix $n$. Order $V(T_n)$ according to increasing distance from $o$ (breaking ties in an arbitrary
deterministic manner). Think of the measure $\mu$ on $V(T_n)$ as being defined sequentially
starting from $o$. $o$ is labelled 1 with probability $\pi(1) \ge \rho$ and then all later vertices
are labelled 1 with probability either $P(0,1)$ or $P(1,1)$ depending on 
the state of the vertice's unique neighbor which has already been assigned a state.
Since both of these probabilities are at least $\rho$, 
we can couple sequentially as in Theorem \ref{thm:2.14} and conclude that $\mu$ restricted to
$V(T_n)$ dominates $\nu_\rho$. As $n$ is arbitrary, we are done.
For (2) implies (3), one observes that by definition,
$$
\mu\{\eta\equiv 0\text{ on } V(T_n)\}= \pi(0) P(0,0)^{|V(T_n)|-1} \,\, \forall n.
$$
(2) now immediately yields that $P(0,0) \le 1-\rho$ which is 
simply (3). $\blacksquare$

\medskip
\noindent
Before giving the proof of Proposition \ref{prop:tree}, we need to summarize some
facts all of which are in (Georgii (1988), pages 247-255). For $J\ge 0$, we let
$$
f_J(t)=\log\big[\frac{\cosh(J+t)} {\cosh(t-J)}\big]
$$
map $R$ to $R$. $f_J$ is
an odd function and concave on $[0,\infty)$. 0 is the unique fixed point
if and only if there is a unique Gibbs state. Otherwise, the fixed points are
$0, t_J$ and $-t_J$ with $t_J>0$. Next, the plus measure, the minus measure and
the free are all tree-indexed Markov chains as defined earlier. (These are all
distinct if there is more than one Gibbs state; this is not always true
in the presence of an external field.) Their respective
transition matrices, denoted by $P^{+,J}$, $P^{-,J}$, and $P^{f,J}$ are given by
$$
\left(
\begin{array}{cc}
P^{+,J}(0,0) & P^{+,J}(0,1) \\
P^{+,J}(1,0) & P^{+,J}(1,1)
\end{array} \right)  =
\left(
\begin{array}{cc}
\frac{e^{J-t_J}} {2\cosh(J-t_J)} & \frac{e^{t_J-J}}{2\cosh(J-t_J)} \\
\frac{e^{-J-t_J}} {2\cosh(J+t_J)} & \frac{e^{J+t_J}}{2\cosh(J+t_J)}
\end{array} \right)
$$

$$
\left(
\begin{array}{cc}
P^{-,J}(0,0) &P^{-,J}(0,1) \\
P^{-,J}(1,0) & P^{-,J}(1,1)
\end{array} \right)  =
\left(
\begin{array}{cc}
\frac{e^{J+t_J}} {2\cosh(J+t_J)} & \frac{e^{-t_J-J}}{2\cosh(J+t_J)} \\
\frac{e^{-J+t_J}} {2\cosh(J-t_J)} & \frac{e^{J-t_J}} {2\cosh(J-t_J)}
\end{array} \right) 
$$

$$
\left(
\begin{array}{cc}
P^{f,J}(0,0) & P^{f,J}(0,1) \\
P^{f,J}(1,0) & P^{f,J}(1,1)
\end{array} \right)  =
\left(
\begin{array}{cc}
\frac{e^{J}} {2\cosh(J)} & \frac{e^{-J}}{2\cosh(J)} \\
\frac{e^{-J}} {2\cosh(J)} & \frac{e^{J}} {2\cosh(J)}
\end{array} \right)
$$

\noindent
{\sl Proof of Proposition \ref{prop:tree}:}  \\
Looking at the formulas for $P^{+,J}$, $P^{-,J}$ and $P^{f,J}$ given above,
one sees that \\
$P^{+,J}(0,1)\le P^{+,J}(1,1)$,
$P^{-,J}(0,1)\le P^{-,J}(1,1)$ and $P^{f,J}(0,1)\le P^{f,J}(1,1)$. 
If $\mu^{J,+}\neq \mu^{J,-}$, then $t_J >0$ and therefore one sees (by looking
at the matrices) that
$$
P^{+,J}(0,1) > P^{f,J}(0,1) >P^{-,J}(0,1).
$$
The result now follows from Proposition \ref{prop:treeFKG}. $\blacksquare$

\noindent{\bf Remarks.} (1) To show only that $\mu^{J,+}$ and $\mu^{J,-}$
dominate different product measures, one does not need the explicit form of the above matrices
but rather only the fact that these measures are tree indexed Markov chains,
Proposition \ref{prop:treeFKG} and an elementary symmetry argument. \\
(2) Using Proposition \ref{prop:treeFKG} and the above form
of the matrices, we immediately see the optimal product measures which
the plus, minus and free measures dominate.

Before proving Proposition \ref{prop:differentJ}, we need the following two lemmas.

\begin{lemma}   \label{lemma:holley}
Given two transition matrices on two states, $P$ and $Q$, let 
$\mu_P$ and $\mu_Q$ 
be the corresponding tree indexed Markov chains on $T$. Then 
$\mu_P$ dominates $\mu_Q$ iff $P(0,1)\ge Q(0,1)$ and
$P(1,1)\ge Q(1,1)$.
\end{lemma}

\noindent
{\sl Proof:} The ``if'' direction is analogous to (3) implies (1) in Theorem
\ref{thm:FKGd} and is just done by coupling sequentially.
For the ``only if'' part,
let $T_n$ be as in Proposition \ref{prop:treeFKG} and observe that
$$
\mu_P\{\eta\equiv 0\text{ on } V(T_n)\}= \pi(0) P(0,0)^{|V(T_n)|-1} \,\, \forall n.
$$
and
$$
\mu_Q\{\eta\equiv 0\text{ on } V(T_n)\}= \pi(0) Q(0,0)^{|V(T_n)|-1} \,\, \forall n.
$$
$\mu_P$ dominating $\mu_Q$ therefore yields $P(0,0)\le Q(0,0)$ or equivalently
$P(0,1)\ge Q(0,1)$. Similarly, by looking at the event 
$\{\eta\equiv 1\text{ on } V(T_n)\}$, one shows that $P(1,1)\ge Q(1,1)$.
$\blacksquare$

\begin{lemma}   \label{lemma:fixedpoint}
If $J_c\le J_1< J_2$, then $t_{J_2}-J_2 \ge t_{J_1}-J_1$. ($t_{J_c}$ is defined to be 0.)
\end{lemma}

\noindent
{\sl Proof:} 
Write $f(J,t)$ for $f_J(t)$, $t(J)$ for $t_J$ and
use subscripts to denote partial derivatives.
Then
$$f_1(J,t)=\tanh(J+t)-\tanh(J-t),$$
and
$$f_2(J,t)=\tanh(J+t)+\tanh(J-t).$$
Differentiate the relation
$$
f(J,t(J))=t(J)
$$
with respect to $J$ and solve to get
$$
t'(J)=\frac{f_1(J,t(J))}{1-f_2(J,t(J))}.
$$
To get $t'(J)\geq 1$, we need:
$$
f_2(J,t(J))<1\quad\text{and}\quad f_1(J,t(J))+f_2(J,t(J))\geq 1.
$$
The first statement is immediate because as a function of
$t$, $f$ crosses the line $y=x$ from above to below. For
the second, note that
$$
f_1(J,t)+f_2(J,t)=2\tanh(J+t)
$$
and hence $f_1+f_2$ is increasing in both variables.
However, since $f_2(J,0)=2\tanh(J),$ $J_c$ is determined
by $\tanh(J_c)=1/2$ and hence it follows that
$$
f_1(J_c,0)+f_2(J_c,0)=1.
$$
Since $f_1+f_2$ is increasing in both variables, we obtain
$$
f_1(J,t)+f_2(J,t)\geq 1
$$
for $J\geq J_c$ and $t\geq 0$.
$\blacksquare$

\noindent
{\sl Proof of Proposition \ref{prop:differentJ}:}  \\
For (i), we have, using Lemma \ref{lemma:fixedpoint} and the exact form of our
matrices,
$$
\frac{P^{+,J_2}(1,1)}{P^{+,J_2}(1,0)}=e^{2t_{J_2}+2J_2} \ge
e^{2t_{J_1}+2J_1} = \frac{P^{+,J_1}(1,1)}{P^{+,J_1}(1,0)}
$$
and
$$
\frac{P^{+,J_2}(0,1)}{P^{+,J_2}(0,0)}=e^{2t_{J_2}-2J_2} \ge
e^{2t_{J_1}-2J_1} = \frac{P^{+,J_1}(0,1)}{P^{+,J_1}(0,0)}
$$
which give
$$
P^{+,J_2}(1,1) \ge P^{+,J_1}(1,1), \text{ and }
P^{+,J_2}(0,1) \ge P^{+,J_1}(0,1).
$$
Now apply Lemma \ref{lemma:holley}.
For (ii), arguing exactly as above with Lemma \ref{lemma:holley},
one can check that for $0\le J_1 \le J_c < J_2$, 
$\mu^{J_2,+}$ dominates $\mu^{J_1,+}$ if and only if $J_1\ge J_2-t_{J_2}$.
This shows the first part of (ii) with $\alpha(J_2)=J_2-t_{J_2}$. 
It then follows that the smallest $J_2 > J_c$ satisfying
$\alpha(J_2)=0$ corresponds exactly to that value of $J >0$ such that
$t_J=J$. A simple computation, using the fixed point equation, shows that
this value is precisely the unique real root of $x^3-x^2-x-1$.
For (iii), one easily checks that for $J$ large,
$f_J(1.1 J) > 1.1J$ and hence for such $J$ that $t_J > 1.1J$.
The result now easily follows from
Lemma \ref{lemma:holley} together with the explicit forms of the matrices
$P^{+,J}$.
$\blacksquare$

\noindent{\bf Remark.} Many of the results presented in this section can be extended 
to the case when there is a nonzero external field.

\section{Exchangeability, FKG and Domination}\label{sec:six}
In this section, we prove Theorem \ref{thm:exchange}.

\noindent
{\sl Proof of Theorem \ref{thm:exchange}:} 
Letting
$$
u_i=P(\eta_1=\cdots=\eta_i=1,\eta_{i+1}=\cdots=\eta_n=0),\quad 0\leq i\leq n,
$$
the FKG lattice condition (see page 78 of Liggett (1985)) then becomes
\begin{equation}\label{eqn:FKG}
u_i^2\leq u_{i-1}u_{i+1},\quad 0<i<n.
\end{equation}
Now, one direction is immediate. If the distribution of
$(\eta_1,...,\eta_n)$ dominates $\nu_{\rho}$, then
$$u_0=P(\eta_1=\eta_2=\cdots=\eta_n=0)\leq (1-\rho)^n.$$
The other direction is harder. We need to prove that if (\ref{eqn:FKG}) and (\ref{4.2new}) hold, 
then 
\begin{equation}\label{eqn:S6tag3}
Eh(\eta_1,...,\eta_n)\geq\int hd\nu_{\rho}
\end{equation}
for all increasing functions $h$ on $\{0,1\}^n$. Since both distributions
appearing in (\ref{eqn:S6tag3}) are exchangeable, it is enough to prove (\ref{eqn:S6tag3}) for 
symmetric increasing functions $h$. To see that this is enough, we
need to check that if $h$ is an increasing function on $\{0,1\}^n$, then
so is its symmetrization $h^*$. Letting $|\eta|=\eta_1+\cdots+\eta_n$, this
is defined by
$$h^*(\eta)=\sum_{\zeta:|\zeta|=k}h(\zeta)\bigg/\binom nk,\quad |\eta|=k.$$
The monotonicity of $h^*$ is equivalent to
$$(n-k)\sum_{\zeta:|\zeta|=k}h(\zeta)\leq(k+1)\sum_{\zeta:|\zeta|=k+1}h(\zeta).$$
To check this inequality, define $\eta^i$ to be the element of $\{0,1\}^n$ obtained
from $\eta$ by flipping the $i$th coordinate. If $h$ is increasing, then
for $\eta$ satisfying $|\eta|=k$,
$$(n-k)h(\eta)\leq\sum_{i:\eta_i=0}h(\eta^i).$$
Summing over all such $\eta$ and changing the order of summation gives
$$(n-k)\sum_{\eta:|\eta|=k}h(\eta)\leq \sum_{i=1}^n\sum_{\eta:|\eta|=k,\eta_i=0}
h(\eta^i)=\sum_{i=1}^n\sum_{\zeta:|\zeta|=k+1,\zeta_i=1}h(\zeta)=(k+1)
\sum_{\zeta:|\zeta|=k+1}h(\zeta)$$
as required.

So, we need to prove (\ref{eqn:S6tag3}) for functions $h$ of the form
$$h(\eta_1,...,\eta_n)=H(\eta_1+\cdots+\eta_n),$$
where $H$ is an increasing function on $\{0,...,n\}$. For this, it
is enough to take $H$ of the form
$$
 H(i)=\left\{ \begin{array}{ll}
1&\text{if }i\geq k\\
0&\text{if }i<k
\end{array} 
\right.
$$
for some $k$. Thus, we need to prove that
\begin{equation}\label{eqn:S6tag4}
\sum_{i=k}^n\binom niu_i\geq \sum_{i=k}^n\binom
ni\rho^i(1-\rho)^{n-i},
\end{equation}
since the left side above is $P(\eta_1+\cdots+\eta_n\geq k)$.

Now, write (\ref{eqn:FKG}) in the form
$$u_i^2\leq\bigg[\frac {\rho}{1-\rho}u_{i-1}\bigg]\bigg[\frac {1-\rho}{\rho}u_{i+1}\bigg]$$
and use the arithmetic-geometric mean inequality to
get
$$2u_i\leq \frac {\rho}{1-\rho}u_{i-1}+\frac {1-\rho}{\rho}u_{i+1}.$$
Dividing by $\rho^i(1-\rho)^{n-i}$
gives
\begin{equation}\label{eqn:S6tag5}
2v_i\leq v_{i-1}+v_{i+1},
\end{equation}
where 
$$
v_i=\frac{u_i}{\rho^i(1-\rho)^{n-i}}.
$$
In other words, the sequence $v_i$ is convex.

We will prove (\ref{eqn:S6tag4}) by contradiction. Suppose it fails for some $k$. Then
for that $k$,
$$\sum_{i=k}^n\binom ni\rho^i(1-\rho)^{n-i}v_i< \sum_{i=k}^n\binom
ni\rho^i(1-\rho)^{n-i}.$$ 
It follows that for some $j\geq k$, $v_j\leq 1$. By (\ref{4.2new}) (which says that
$v_0\leq 1$) and (\ref{eqn:S6tag5}), $v_i\leq 1$ for all $0\leq i\leq k$. This gives
$$\sum_{i=0}^{k-1}\binom ni\rho^i(1-\rho)^{n-i}v_i\leq \sum_{i=0}^{k-1}\binom
ni\rho^i(1-\rho)^{n-i}.$$ 
Adding the last two displays yields
$$\sum_{i=0}^n\binom niu_i=
\sum_{i=0}^n\binom ni\rho^i(1-\rho)^{n-i}v_i<
\sum_{i=0}^n\binom ni\rho^i(1-\rho)^{n-i}.$$ 
But this is a contradiction, since the two extreme sums above are 
equal to 1.$\blacksquare$ 

\noindent {\bf Remark.} The statement of Theorem 1.3 is not true if either
assumption of exchangeability or FKG is omitted, even for $n=2$. For counterexamples,
suppose first that exchangeability is omitted. Then one can take $\eta_1,\eta_2$
to be independent with $P(\eta_1=1)=\alpha$ and $P(\eta_2=1)=\beta$. The FKG
condition holds for all $\alpha$ and $\beta$. However, the distribution of $(\eta_1,\eta_2)$
stochastically dominates $\nu_{\rho}$ if and only if $\min (\alpha,\beta)\geq\rho$, while 
(\ref{4.2new}) holds if and only if $(1-\alpha)(1-\beta)\leq(1-\rho)^2$. Suppose now that the FKG
assumption is omitted, and take $P(\eta_1=1,\eta_2=0)=P(\eta_1=0,\eta_2=1)=1/2.$
This is exchangeable and its distribution does not dominate any nontrivial
product measure, yet (\ref{4.2new}) is satisfied for all $\rho$.

\section{An example} \label{sec:eight}

In this section, we present a example of $X_1,X_2,X_3,X_4$ which are 0,1 valued,
exchangeable and FKG but which are not extendible to an infinite exchangeable process.
As usual, let $u_i$ be the probability of a
configuration with $i$ ones and $4-i$ zeros. Take
$$
u_0=u_4=c\lambda^2,\quad u_1=u_3=c\lambda,\quad u_2=c,
$$
where
$$
c=\frac 1{2\lambda^2+8\lambda+6}.
$$
This satisfies the FKG condition iff $\lambda\geq 1$. If the measure
were infinitely extendible, there would be a random variable
$0\leq W\leq 1$ so that
$$
u_i=EW^i(1-W)^{4-i}.
$$
Then
$$
E[W^2-\lambda W(1-W)]^2=0,
$$
so that $W$ can take on only the values $0$ and $\lambda/(1+\lambda)$. 
Similarly, $W$ can take on only the values $1$ and $1/(1+\lambda)$.
This is a contradiction unless $\lambda=1$.

\section{Some open questions} \label{sec:nine}

The first five questions concern the contact process and the following
two questions concern the Ising model.

\noindent
1. Fix $d\ge 1$. Given $\rho >0$ does there exist $\lambda >\lambda_c$ such that 
$\nu_{\rho}$ stochastically dominates $\overline\nu_{\lambda}$?
(This would be an essential strengthening of the fact that the critical
contact process dies out.)

\noindent
2. For $d\ge 2$, does there exist $\lambda > \lambda_c$ such that
$\overline\nu_{\lambda}$ does not percolate?
(In words, is the critical value for percolation different than the usual critical value?)

\medskip
\noindent
Observe that a positive answer to question 1 would yield a positive answer to question 2.

\noindent
3. For bounded degree graphs $G$ with site percolation critical value less than 1, 
does there exist $\lambda$ such that for the symmetric contact process on $G$ with parameter 
$\lambda$, $\overline\nu_{\lambda}$ percolates?

\noindent
4. For bounded degree graphs, is it the case that for all $\rho< 1$,
there exists $\lambda$ such that for the symmetric contact process on $G$ with parameter 
$\lambda$, $\overline\nu_{\lambda}$ stochastically dominates $\nu_\rho$?

\noindent
5. Assume that for the parameter $\lambda$, 
$\overline\nu_{\lambda}$ stochastically dominates $\nu_\rho$
for the symmetric contact process on $Z^+$. Does it follow that
for any bounded degree graph $G$, $\overline\nu_{\lambda}$
for the corresponding symmetric contact process on $G$ also
dominates $\nu_\rho$?

\medskip
\noindent
Observe that a positive answer to question 5 implies a 
positive answer to question 4 which in turn implies a 
positive answer to question 3.

\medskip
\noindent
{\bf Remark.}  An interesting test case for question 5 is $Z^+$ with $n$ dangling edges;
that is, the example studied in Proposition \ref{prop:3.1}. We have seen that question 4 holds in 
this case for any $n$ and we have seen that question 5 fails for an asymmetric version.

\medskip
\noindent
6. Given any nonamenable transitive graph, does the plus state for the Ising model
for large values of $J$ dominate high density product measures?

\noindent
7. Is amenability for transitive graphs characterized by the property that the plus and minus 
states for the Ising model for fixed $J$ dominate the same set of product measures or alternatively 
by the property that the plus states for different $J$'s cannot be stochastically ordered?

\noindent
8. Is there some reasonable version of 
Theorem \ref{thm:FKGd} and Proposition \ref{prop:treeFKG}
for Markov fields on trees which are not tree-indexed Markov chains?

\noindent
{\bf Acknowledgement:}  We thank the referee for an extremely careful
reading and for a number of suggestions.

\end{document}